\documentclass{amsart}
\usepackage{setspace}

\usepackage{amsmath,amsthm,amsfonts,amscd,amssymb}

\newtheorem{thm}{Theorem}[section]
\newtheorem{cor}[thm]{Corollary}
\newtheorem{lem}[thm]{Lemma}

\theoremstyle{definition}

\theoremstyle{remark}
\newtheorem{rem}{Remark}[section]

\begin{document}

\title{Effective structure theorems for symplectic spaces via height}
\author{Lenny Fukshansky}

\address{Department of Mathematics, 850 Columbia Avenue, Claremont McKenna College, CA 91711}
\email{lenny@cmc.edu}

\begin{abstract}
Given a $2k$-dimensional symplectic space $(Z,F)$ in $N$ variables, $1 < 2k \leq N$, over a global field  $K$, we prove the existence of a symplectic basis for $(Z,F)$ of bounded height. This can be viewed as a version of Siegel's lemma for a symplectic space. As corollaries of our main result, we prove the existence of a small-height decomposition of $(Z,F)$ into hyperbolic planes, as well as the existence of two generating flags of totally isotropic subspaces. These present analogues of known results for quadratic spaces. A distinctive feature of our argument is that it works simultaneously for essentially any field with a product formula, algebraically closed or not. In fact, we prove an even more general version of these statements, where canonical height is replaced with twisted height. All bounds on height are explicit.
\end{abstract}

\maketitle

\def\A{{\mathcal A}}
\def\B{{\mathcal B}}
\def\C{{\mathcal C}}
\def\D{{\mathcal D}}
\def\F{{\mathcal F}}
\def\x{{\mathcal H}}
\def\I{{\mathcal I}}
\def\J{{\mathcal J}}
\def\K{{\mathcal K}}
\def\kk{{\mathfrak K}}
\def\L{{\mathcal L}}
\def\M{{\mathcal M}}
\def\mm{{\mathfrak m}}
\def\MM{{\mathfrak M}}
\def\R{{\mathcal R}}
\def\s{{\mathcal S}}
\def\V{{\mathcal V}}
\def\X{{\mathcal X}}
\def\Y{{\mathcal Y}}
\def\Z{{\mathcal Z}}
\def\H{{\mathcal H}}
\def\OO{{\mathfrak O}}
\def\CC{{\mathfrak C}}
\def\cee{{\mathbb C}}
\def\Nn{{\mathbb N}}
\def\pee{{\mathbb P}}
\def\que{{\mathbb Q}}
\def\real{{\mathbb R}}
\def\zed{{\mathbb Z}}
\def\hyp{{\mathbb H}}
\def\AA{{\mathbb A}}
\def\hs{{\hat \sigma}}
\def\gmn{{\mathbb G_m^N}}
\def\qbar{{\overline{\mathbb Q}}}
\def\kkbar{{\overline{\mathfrak K}}}
\def\kbar{{\overline{K}}}
\def\eps{{\varepsilon}}
\def\vek{{\varepsilon_k}}
\def\ahat{{\hat \alpha}}
\def\bhat{{\hat \beta}}
\def\gt{{\tilde \gamma}}
\def\h{{\tfrac12}}
\def\ba{{\boldsymbol a}}
\def\be{{\boldsymbol e}}
\def\bei{{\boldsymbol e_i}}
\def\bc{{\boldsymbol c}}
\def\bm{{\boldsymbol m}}
\def\bk{{\boldsymbol k}}
\def\bi{{\boldsymbol i}}
\def\bl{{\boldsymbol l}}
\def\bq{{\boldsymbol q}}
\def\bu{{\boldsymbol u}}
\def\bt{{\boldsymbol t}}
\def\bs{{\boldsymbol s}}
\def\bv{{\boldsymbol v}}
\def\bw{{\boldsymbol w}}
\def\bx{{\boldsymbol x}}
\def\bX{{\boldsymbol X}}
\def\bz{{\boldsymbol z}}
\def\bwy{{\boldsymbol y}}
\def\bg{{\boldsymbol g}}
\def\bY{{\boldsymbol Y}}
\def\bL{{\boldsymbol L}}
\def\baa{{\boldsymbol\alpha}}
\def\bb{{\boldsymbol\beta}}
\def\bet{{\boldsymbol\eta}}
\def\bxi{{\boldsymbol\xi}}
\def\bo{{\boldsymbol 0}}
\def\bol{{\boldsymbol 1}_L}
\def\ep{\varepsilon}
\def\p{\boldsymbol\varphi}
\def\q{\boldsymbol\psi}
\def\rank{\operatorname{rank}}
\def\aut{\operatorname{Aut}}
\def\lcm{\operatorname{lcm}}
\def\sgn{\operatorname{sgn}}
\def\spn{\operatorname{span}}
\def\md{\operatorname{mod}}
\def\Norm{\operatorname{Norm}}
\def\dim{\operatorname{dim}}
\def\det{\operatorname{det}}
\def\Vol{\operatorname{Vol}}
\def\rk{\operatorname{rk}}
\def\GL{\operatorname{GL}}

\section{Introduction}
Throughout this paper, we let $K$ be either a number field, a function field (i.e. a finite algebraic extension of the field of rational functions in one variable over an arbitrary field), or the algebraic closure of one or the other. Let
\begin{equation}
\label{form}
F(\bX,\bY) = \sum_{i=1}^N \sum_{j=1}^N f_{ij} X_i Y_j
\end{equation}
be an alternating bilinear form in $N \geq 2$ variables with coefficients in $K$. We will also write $F=(f_{ij})_{1 \leq i,j \leq N}$ for the anti-symmetric $N \times N$ coefficient matrix of $F$, i.e. $f_{ij}=-f_{ji}$ for all $1 \leq i,j \leq N$. In particular, $f_{ii}=0$ for all $1 \leq i \leq N$, and the associated quadratic form $F(\bX) = F(\bX,\bX)$ is identically zero on $K^N$. 

Let $Z$ be a $2k$-dimensional subspace of $K^N$, $1 \leq k \leq N/2$, and let us write $(Z,F)$ for the symplectic space defined on $Z$ by $F$. We will assume that $(Z,F)$ is regular, meaning that for every $\bo \neq \bx \in Z$ there exists $\bwy \in Z$ such that $F(\bx,\bwy) \neq 0$. Then $(Z,F)$ has a symplectic basis (see for instance \cite{scharlau}), that is a basis $\bx_1,\dots,\bx_k,\bwy_1,\dots,\bwy_k$ for $Z$ over $K$ such that
\begin{equation}
\label{sympl_basis}
F(\bx_i,\bx_j) = F(\bwy_i, \bwy_j) = F(\bx_i,\bwy_j) = 0\ \forall\ 1 \leq i \neq j \leq k,\ F(\bx_i,\bwy_i) = 1\ \forall\ 1 \leq i \leq k.
\end{equation}
A subspace $V$ of $Z$ is called totally isotropic if $F(\bx,\bwy)=0$ for all $\bx,\bwy \in V$, and a maximal totally isotropic subspace of $(Z,F)$ is called a Lagrangian. All Lagrangians of $(Z,F)$ have the same dimension; it is an easy consequence of (\ref{sympl_basis}) that this dimension is $k$. Indeed, it is easy to see that $V_1=\spn_K \{ \bx_1,\dots,\bx_k\}$ and $V_2=\spn_K \{ \bwy_1,\dots,\bwy_k\}$ are Lagrangians in $(Z,F)$. Moreover, $(Z,F)$ is a hyperbolic space over these Lagrangians, meaning that
\begin{equation}
\label{hyperbolic}
Z = \hyp_1 \perp \dots \perp \hyp_k,
\end{equation}
where for each $1 \leq i \leq k$, $\hyp_i = \spn_K \{ \bx_i,\bwy_i \}$ is a hyperbolic plane, $\perp$ stands for orthogonal direct sum, and orthogonality throughout this paper is always meant with respect to $F$. This means that once we know how to find a symplectic basis for $(Z,F)$, we immediately obtain two Lagrangians as well as an orthogonal decomposition of $(Z,F)$ into hyperbolic planes. However the classical result about the existence of a basis satisfying (\ref{sympl_basis}) is ineffective, i.e. it provides no information as to how does one find such a basis.

The main goal of this paper is to prove an effective version of the existence theorem for a symplectic basis, and derive from it effective statements about existence of Lagrangians and a hyperbolic decomposition for a regular symplectic space. We use the approach of height functions, which will be formally introduced in section~2. We will define a height function $H$ on the points of a projective space over $K$, and in particular will talk about height of vectors and subspaces of $K^N$ to mean height of the corresponding projective points; specifically, subspaces of $K^N$ will be viewed as points on a corresponding Grassmanian. We will also give a slightly different definition for the height $\H$ of our alternating bilinear form $F$. Loosely speaking, height measures the arithmetic complexity of objects in question, meaning that the smaller is the height of a projective point the less ``arithmetically complicated'' this point is. In particular, height satisfies the crucial finiteness property: any set of projective algebraic points of bounded height and degree is always finite (this will be rigorously discussed in section 2, especially see (\ref{northcott})). Therefore, proving the existence of a point or subspace of bounded height over $K$ that satisfies some arithmetic conditions may provide a search bound for points satisfying such conditions. Hence our goal will be to prove effective theorems for symplectic spaces in the sense of providing bounds on height. We can now state our main result.

\begin{thm} \label{main} Let $(Z,F)$ be a regular $2k$-dimensional symplectic space in $N$ variables over $K$, where $1 \leq k < 2k \leq N$. Then there exists a symplectic basis $\bx_1,\dots,\bx_k,\bwy_1,\dots,\bwy_k$ for $Z$ satisfying (\ref{sympl_basis}) such that
\begin{equation}
\label{effective_basis}
\prod_{i=1}^{2k} H(\bx_i)H(\bwy_i) \leq \left( C_K(N,2k) H(Z) \right)^{a_k} \H(F)^{b_k},
\end{equation}
where $C_K(N,2k)$ is a field constant defined in section~2 below,
\[ a_k = \left\{ \begin{array}{ll}
\frac{k^2+4k}{4} & \mbox{if $2|k$} \\
\frac{k^2+4k-1}{4} & \mbox{if $2 \nmid k$},
\end{array}
\right. \]
and
\[ b_k = \left\{ \begin{array}{ll}
\frac{2k^3+9k^2-14k}{12} & \mbox{if $2|k$} \\
\frac{2k^3+9k^2-14k+3}{12} & \mbox{if $2 \nmid k$}.
\end{array}
\right. \]
\end{thm}

\noindent
An immediate corollary of Theorem \ref{main} is an effective version of Witt decomposition for $(Z,F)$, which in the symplectic case is just a decomposition into hyperbolic planes.

\begin{cor} \label{hyp_dec} Let the notation be as in Theorem \ref{main}, then there exists a decomposition (\ref{hyperbolic}) for $(Z,F)$ with
\begin{equation}
\label{effective_hyp}
\prod_{i=1}^{2k} H(\hyp_i) \leq \left( C_K(N,2k) H(Z) \right)^{a_k} \H(F)^{b_k}.
\end{equation}
\end{cor}

\proof
For each $1 \leq i \leq k$, take $\hyp_i = \spn_K \{ \bx_i,\bwy_i \}$, then by Lemma \ref{wedge} below $H(\hyp_i) \leq H(\bx_i) H(\bwy_i)$, and the statement of the corollary follows from (\ref{effective_basis}).
\endproof

\noindent
We can now also establish the existence of flags of totally isotropic subspaces of bounded height, whose union generates~$Z$.

\begin{cor} \label{flags} Let the notation be as in Theorem \ref{main}. For each $1 \leq n \leq k$, there exist totally isotropic subspaces $V_n$ and $W_n$ of $(Z,F)$ such that $\dim_K V_n = \dim_K W_n = n$, $V_n \cap W_n = \{ \bo \}$,
\begin{equation}
\label{flag}
V_1 \subset V_2 \subset \dots \subset V_k,\ W_1 \subset W_2 \subset \dots \subset W_k,
\end{equation}
and
\begin{equation}
\label{flag_bound}
H(V_n) H(W_n) \leq \left( C_K(N,2k)^{a_k} H(Z)^{a_k} \H(F)^{b_k} \right)^{\frac{n}{k}}.
\end{equation}
In particular, $(Z,F)$ is generated by the two small-height Lagrangians $V_k$ and $W_k$, i.e. $Z = \spn_K \{ V_k, W_k \}$.
\end{cor}

\proof
With notation of Theorem \ref{main}, assume without loss of generality that the symplectic basis vectors are ordered in such a way that
$$H(\bx_1) H(\bwy_1) \leq H(\bx_2) H(\bwy_2) \leq \dots H(\bx_k) H(\bwy_k).$$
Then let $V_n = \spn_K \{ \bx_1, \dots, \bx_n \}$ and $W_n = \spn_K \{ \bwy_1, \dots, \bwy_n \}$, for each $1 \leq n \leq k$, and notice that by Lemma \ref{wedge},
$$H(V_n) H(W_n) \leq H(\bx_1) \dots H(\bx_n) H(\bwy_1) \dots H(\bwy_n).$$
The statement of the corollary now follows from (\ref{effective_basis}).
\endproof

These results should be viewed as symplectic space analogues of Siegel's lemma and effective decomposition theorems for quadratic spaces. The name Siegel's lemma usually refers to results about  the existence of a basis of small height for a vector space over a global field (see \cite{vaaler:siegel}, \cite{thunder}, and \cite{absolute:siegel}). In case of a quadratic space (i.e. when $F$ is a symmetric bilinear form), a version of Siegel's lemma with additional conditions, asserting the existence of an orthogonal basis of small height, has been proved in \cite{me:witt} over a number field and in \cite{me:qbar_quad} over $\qbar$. Theorem \ref{main} is precisely a symplectic space analogue of these theorems. 

There has been a large number of results on small-height zeros of quadratic forms, starting with a classical theorem of Cassels \cite{cassels:small}. One of the directions generalizing Cassels' theorem produced results on small-height linear subspaces of a quadratic space on which the quadratic form vanishes identically (see \cite{schlickewei}, \cite{schmidt:schlickewei}, \cite{vaaler:smallzeros}, and \cite{vaaler:smallzeros2}). Corollary \ref{flags} should be viewed as an analogue of these results for a symplectic space. Finally, the structural results for a quadratic space, such as effective Witt decomposition, have been proved in \cite{me:witt} and \cite{me:qbar_quad}; Corollary \ref{hyp_dec} serves as a symplectic space analogue of this.

In the case of a quadratic space, such problems were usually treated separately by different methods  over a number field, function field, or algebraic closures. The distinctive feature of the symplectic situation is that, because it is much more linear, we are able to treat these problems at once over any global field with a product formula for which a Siegel's lemma type result exists - this is due to the purely combinatorial nature of our argument. Moreover, we prove our main result in terms of more general twisted heights (see Theorem \ref{main_twisted}), from which Theorem \ref{main} follows immediately.

This paper is structured as follows. In section~2 we set the notation and define the height functions, and present a few technical lemmas on properties of heights. In section~3 we prove a combinatorial lemma (Lemma \ref{combin_lemma}), which we later use to obtain Theorem \ref{main}. In section~4 we derive Theorem \ref{main} by means of proving the more general Theorem \ref{main_twisted} with the use of Siegel's lemma, stated as Theorem \ref{siegel}, and Lemma \ref{combin_lemma}.

\section{Notation and heights}

We start with some notation, following \cite{absolute:siegel}. Throughout this paper, $K$ will either be a number field (finite extension of $\que$), a function field, or algebraic closure of one or the other; in fact, for the rest of this section, unless explicitly specified otherwise, we will assume that $K$ is either a number field of a function field, and will write $\kbar$ for its algebraic closure. By a function field we will always mean a finite algebraic extension of the field $\kk = \kk_0(t)$ of rational functions in one variable over a field $\kk_0$, where $\kk_0$ can be any field. When $K$ is a number field, clearly $K \subset \kbar = \qbar$; when $K$ is a function field, $K \subset \kbar = \kkbar$, the algebraic closure of $\kk$. In the number field case, we write $d = [K:\que]$ for the global degree of $K$ over $\que$; in the function field case, the global degree is $d = [K:\kk]$, and we also define the effective degree of $K$ over $\kk$ to be
$$\mm(K,\kk) = \frac{[K:\kk]}{[K_0:\kk_0]},$$
where $K_0$ is the algebraic closure of $\kk_0$ in $K$. If $K$ is a number field, we let $\D_K$ be its discriminant; if $K$ is a function field, we will also write $g(K)$ for the genus of $K$, as defined by the Riemann-Roch theorem (see \cite{thunder} for details). We can now define the field constant $C_K(N,L)$, which appears in our upper bounds:
\[ C_K(N,L) = \left\{ \begin{array}{ll}
N^{\frac{L}{2}} |\D_K|^{\frac{L}{2d}} & \mbox{if $K$ is a number field} \\
\exp\left(\frac{g(K)-1+\mm(K,\kk)}{\mm(K,\kk)}\right) & \mbox{if $K$ is a function field} \\
3^{\frac{L(L-1)}{2}} & \mbox{if $K = \qbar$} \\
2 & \mbox{if $K = \overline{\kk}$},
\end{array}
\right. \]

Next we discuss absolute values on $K$. Let $M(K)$ be the set of places of $K$. For each place $v \in M(K)$ we write $K_v$ for the completion of $K$ at $v$ and let $d_v$ be the local degree of $K$ at $v$, which is $[K_v:\que_v]$ in the number field case, and $[K_v:\kk_v]$ in the function field case. 

If $K$ is a number field, then for each place $v \in M(K)$ we define the absolute value $|\ |_v$ to be the unique absolute value on $K_v$ that extends either the usual absolute value on $\real$ or $\cee$ if $v | \infty$, or the usual $p$-adic absolute value on $\que_p$ if $v|p$, where $p$ is a rational prime. 

If $K$ is a function field, then all absolute values on $K$ are non-archimedean. For each $v \in M(K)$, let $\OO_v$ be the valuation ring of $v$ in $K_v$ and $\MM_v$ the unique maximal ideal in $\OO_v$. We choose the unique corresponding absolute value $|\ |_v$ such that:
\begin{trivlist}
\item (i) if $1/t \in \MM_v$, then $|t|_v = e$,
\item (ii) if an irreducible polynomial $p(t) \in \MM_v$, then $|p(t)|_v = e^{-\deg(p)}$.
\end{trivlist}

\noindent
In both cases, for each non-zero $a \in K$ the {\it product formula} reads
\begin{equation}
\label{product_formula}
\prod_{v \in M(K)} |a|^{d_v}_v = 1.
\end{equation} 

We can now define local norms on vectors. For each $v \in M(K)$ define a local norm $\|\ \|_v$ on $K_v^N$ by
\[ \|\bx\|_v = \left\{ \begin{array}{ll}
\max_{1 \leq i \leq N} |x_i|_v & \mbox{if $v \nmid \infty$} \\
\left( \sum_{i=1}^N |x_i|_v^2 \right)^{1/2} & \mbox{if $v | \infty$}
\end{array}
\right. \]
for each $\bx \in K_v^N$. We define the following global height function on $K^N$:
\begin{equation}
H(\bx) = \left( \prod_{v \in M(K)} \|\bx\|^{d_v}_v \right)^{1/d},
\end{equation}
for each $\bx \in K^N$.  More generally, let us define the {\it twisted height} on $K^N$ as introduced by J. L. Thunder.  We write $K_{\AA}$ for the ring of adeles of $K$, and view $K$ as a subfield of $K_{\AA}$ under the diagonal embedding (see \cite{weil} for details). Let $A \in GL_N(K_{\AA})$ with local components $A_v \in GL_N(K_v)$. The corresponding twisted height on $K^N$ is defined by
\begin{equation}
H_A(\bx) = \left( \prod_{v \in M(K)} \|A_v \bx\|^{d_v}_v \right)^{1/d},
\end{equation}
for all $\bx \in K^N$. Given any finite extension $E/K$, $K_{\AA}$ can be viewed as a subring of $E_{\AA}$, and let us also write $A$ for the element of $GL_N(E_{\AA})$ which coincides with $A$ on $K_{\AA}^N$. The corresponding twisted height on $E^N$ extends the one on $K^N$, hence $H_A$ is a height on $\kbar$. Notice also that the usual height $H$ as defined above is simply $H_I$, where $I$ is the identity element of $GL_N(K_{\AA})$ all of whose local components are given by $N \times N$ identity matrices. Due to the normalizing exponent $1/d$, our height functions are absolute, i.e. for points over $\qbar$ or $\kkbar$, respectively, their value does not depend on the field of definition. This means that if $\bx$ is in $\qbar^N$ or $\kkbar^N$, then for every $A \in GL_N(K_{\AA})$, $H_A(\bx)$ can be evaluated over any number field or function field, respectively, containing the coordinates of $\bx$, and so $H_A$ provides a height on~$\kbar^N$.

A fundamental property of heights, sometimes referred to as the Northcott property, is that for every $A \in GL_N(K_{\AA})$,
\begin{equation}
\label{northcott}
\left| \left\{ [\baa] \in \pee^{N-1}(\kbar)\ :\ \deg([\baa]) \leq B, H_A(\baa) \leq C \right\} \right| < \infty,
\end{equation}
where $\pee^{N-1}(\kbar)$ is $(N-1)$-dimensional projective space over $\kbar$, $\baa = (\alpha_1,\dots,\alpha_N)$ is in $\kbar^N$, and so $[\baa]$ is the corresponding projective point, $B,C$ are positive real numbers, and $\deg(\baa) = [\que(\alpha_1,\dots,\alpha_N):\que]$ if $K$ is a number field, or $[\kk(\alpha_1,\dots,\alpha_N):\kk]$ if $K$ is a function field, i.e. it is the algebraic degree of  $\baa$ over the ground field over which $K$ is defined. We define $\deg([\baa])$, the degree of the projective point represented by $\baa$, by
$$\deg([\baa]) = \min \{ \deg(\baa') \ : \ \baa' \in \kbar^N,\ [\baa'] = [\baa] \}.$$

We can now extend our notation to define Schmidt twisted height on matrices, which is the same as height function on subspaces of $\kbar^N$. Let $A \in GL_N(K_{\AA})$, $\be_1,\dots,\be_N$ be the standard basis for $K^N$, and $1 \leq J \leq N$. Then $J$-th exterior component $\bigwedge^J K^N$ can be identified with the vector space $K^{\binom{N}{J}}$ via the cannonical isomorphism that sends the wedge products $\be_{i_1} \wedge \dots \wedge \be_{i_J}$, $1 \leq i_1 < \dots < i_J \leq N$, to the standard basis elements of $K^{\binom{N}{J}}$ in lexicographic order. This also identifies $\bigwedge^J A$ with an element of $GL_{\binom{N}{J}}(K_{\AA})$, and so we can talk about the height $H_{\bigwedge^J A}$ on $\bigwedge^J \kbar^N$. Let $X$ be an $N \times J$ matrix of rank $J$ whose column vectors are $\bx_1,...,\bx_J \in K^N$, then we define
$$H_A(X) = H_{\bigwedge^J A}(\bx_1 \wedge\ ...\ \wedge \bx_J).$$
In the same manner, we define the height of a $J \times N$ matrix to be the height of the wedge product of its row vectors instead of column vectors. Now let $V \subseteq \kbar^N$ be a subspace of dimension $J$, $1 \leq J \leq N$, defined over $K$. Choose a basis $\bx_1,...,\bx_J$ for $V$ over $K$, and write $X = (\bx_1 \dots \bx_J)$ for the corresponding $N \times J$ basis matrix. Define the height of $V$ by $H_A(V) = H_A(X)$. This definition is legitimate, since it does not depend on the choice of the basis for $V$: let $\bwy_1,...,\bwy_J$ be another basis for $V$ over $K$ and $Y = (\bwy_ 1 \dots \bwy_J)$ the corresponding $N \times J$ basis matrix, then there exists $W \in GL_J(K)$ such that $Y = XW$, and so
$$\bwy_1 \wedge\ ...\ \wedge \bwy_J = (\det W)\ \bx_1 \wedge\ ...\ \wedge \bx_J,$$
hence, by the product formula 
$$H_A(Y) = H_{\bigwedge^J A}(\bwy_1 \wedge\ ...\ \wedge \bwy_J) = H_{\bigwedge^J A}(\bx_1 \wedge\ ...\ \wedge \bx_J) = H_A(X).$$
On the other hand, there exists an $(N-J) \times N$ matrix $B$ of rank $N-J$ with entries in $K$ such that
\begin{equation}
\label{null_matrix}
V = \left\{ \bx \in \kbar^N : B \bx = 0 \right\}.
\end{equation}
An important duality principle relates heights of $V$ and $B$. For $A \in GL_N(K_{\AA})$ with local components $A_v \in GL_N(K_v)$ for every $v \in M(K)$, let $A^* \in GL_N(K_{\AA})$ be given by the local components $(A_v^t)^{-1} \in GL_N(K_v)$ for every $v \in M(K)$. We also define
$$|\det A|_{\AA} = \left( \prod_{v \in M(K)} |\det A_v|^{d_v}_v \right)^{1/d}.$$
The following is Theorem 1.1 of \cite{absolute:siegel} (see also Duality Theorem in section~2 of \cite{thunder:flag}).

\begin{lem} \label{dual} For any subspace $V \subseteq \kbar^N$ and $A \in GL_N(K_{\AA})$, we have
$$H_{A^*}(B) = |\det A|^{-1}_{\AA} H_A(V),$$
where $B$ is as in (\ref{null_matrix}).
\end{lem}

\noindent
In particular, this implies that $H(V)=H(B)$ if $B$ is as in (\ref{null_matrix}), since clearly for the identity $I \in GL_N(K_{\AA})$, $I^* = I$ and $|\det I|_{\AA} = 1$.

We also define height of our bilinear form $F$ in the following conventional way: let $\H(F)$ be the usual height $H$ of the anti-symmetric matrix $(f_{ij})_{1 \leq i,j \leq N}$, viewed as a vector in $K^{N^2}$. Notice that it is different from the height on matrices defined above, which is why we denote it by $\H$ instead of $H$.

Finally, we define certain dilation constants for an element $A \in  GL_N(K_{\AA})$ that will appear in our bounds (see Lemmas 3.1, 3.2, and Proposition 4.1 of \cite{absolute:siegel}; see also \cite{absolute:siegel1}). Roughly speaking, as we will see in Lemma \ref{comparison} below, these constants indicate by how much does a given automorphism $A$ of $K_{\AA}^N$ "distort" the corresponding twisted height $H_A$ as compared to $H$, the canonical height. Let $A_v = (a^v_{ij})_{1 \leq i,j \leq N} \in  GL_N(K_v)$ be local components of $A$ for each $v \in M(K)$, and let us write $A_v^{-1} = (b^v_{ij})_{1 \leq i,j \leq N}$. Then for all but finitely many places $v \in M(K)$ the corresponding map $A_v$ is an isometry; in fact, let $M_A(K) \subset M(K)$ be the finite (possibly empty) subset of places $v$ at which $A_v$ is {\it not} an isometry. For each $v \notin M_A(K)$, define $\C^v_1(A) = \C^v_2(A) = 1$, and for each $v \in M_A(K)$, let
\begin{equation}
\label{aut_const_loc}
\C^v_1(A) = \left( \sum_{l=1}^N \sum_{m=1}^N |b^v_{lm}|_v \right)^{-1},\ \C_2^v(A) = \sum_{i=1}^N \sum_{j=1}^N |a^v_{ij}|_v.
\end{equation}
Then define
\begin{equation}
\label{aut_const}
\C_1(A) = \prod_{v \in M(K)} (\C^v_1)^{d_v/d},\ \C_2(A) = \prod_{v \in M(K)} (\C^v_2)^{d_v/d},
\end{equation}
both of which are products of only a finite number of non-trivial terms. With this notation, it will also be convenient to define
\begin{equation}
\label{aut_const1}
\CC(A) = \frac{\C_2(A)}{\C_1(A)} = \prod_{v \in M_A(K)} \left( \sum_{i,j,l,m = 1}^N |a^v_{ij} b^v_{lm}|_v \right)^{d_v/d},
\end{equation}
and
\begin{equation}
\label{aut_const1.1}
\CC'(A) = \frac{\CC(A) |\det A|_{\AA}^{1/2}}{\C_1(A)^2}.
\end{equation}
Clearly, in the case when $A = I$ is the identity element of $GL_N(K_{\AA})$, $\CC'(A) = \CC(A) = \C_1(A) = \C_2(A) = 1$. Another important observation is that, since for every $v \in M(K)$, $(A_v^t)^{-1} = (A_v^{-1})^t$, therefore
\begin{eqnarray}
\label{aut_const2}
& & \C^v_1(A^*)^{-1} = \C^v_2(A),\ \C^v_2(A^*) = \C^v_1(A)^{-1},\ \nonumber \\
& & \C_1(A^*)^{-1} = \C_2(A),\ \C_2(A^*) = \C_1(A)^{-1},\ \nonumber \\
& & \CC(A^*) = \CC(A).
\end{eqnarray}
\smallskip

Next we present some technical lemmas that we use later in our main proof, detailing the key properties of height functions. The first one shows that the canonical height $H$ and the twisted height $H_A$ are comparable for each $A \in GL_N(K_{\AA})$ with the comparison constants being precisely the dilation constants $\C_1(A), \C_2(A)$ defined above. This is Proposition 4.1 of \cite{absolute:siegel}. 

\begin{lem} \label{comparison} Let $A \in GL_N(K_{\AA})$. Then
\begin{equation}
\label{comp1}
\C_1(A) H(\bx) \leq H_A(\bx) \leq \C_2(A) H(\bx),
\end{equation}
for all $\bx \in \kbar^N$, where $\C_1(A)$ and $\C_2(A)$ are as in (\ref{aut_const}) above.
\end{lem}

\begin{rem} A simple consequence of Lemma \ref{comparison} and (\ref{aut_const2}) which will be useful to us is that for all $\bx \in \kbar^N$,
\begin{equation}
\label{comp2}
H_{A^*}(\bx) \leq \C_1(A)^{-2} H_A(\bx).
\end{equation}
\end{rem}

\noindent
The next lemma is a consequence of Laplace's expansion, and can be found as Lemma 4.7 of \cite{absolute:siegel} (also see pp. 15-16 of \cite{vaaler:siegel}).

\begin{lem} \label{wedge} Let $X$ be a $N \times J$ matrix over $\kbar$ with column vectors $\bx_1,...,\bx_J$, and let $A \in GL_N(K_{\AA})$. Then
\begin{equation}
\label{prod_1}
H_A(X) = H_{\bigwedge^J A}(\bx_1 \wedge \bx_1\ ...\ \wedge \bx_J) \leq \prod_{i=1}^J H_A(\bx_i).
\end{equation}
More generally, if the $N \times J$ matrix $X$ can be partitioned into blocks as $X = (X_1\ X_2)$, then
 \begin{equation}
\label{prod_2}
H_A(X) \leq H_A(X_1) H_A(X_2).
\end{equation}
\end{lem}
\smallskip

\noindent
The following well known fact is an immediate corollary of Theorem 1 of \cite{vaaler:struppeck} adapted over $\kbar$ and extended to twisted height.

\begin{lem} \label{intersection} Let $U_1$ and $U_2$ be subspaces of $\kbar^N$, and let $A \in GL_N(K_{\AA})$. Then
$$H_A(U_1 \cap U_2) \leq H_A(U_1) H_A(U_2).$$
\end{lem}

\noindent
The next one is a generalization of Lemma 2.3 of \cite{me:witt} over~$\kbar$ and with the twisted height $H_A$ replacing canonical height $H$. We present the proof here for the purposes of self-containment.

\begin{lem} \label{mult} Let $X$ be a $N \times J$ matrix over $\kbar$ with column vectors $\bx_1,...,\bx_J$, $A \in GL_N(K_{\AA})$, and let $F$ be a bilinear form in $N$ variables, as above (we also write $F$ for its $N \times N$ coefficient matrix). Then
\begin{equation}
\label{prod_3}
H_A(F X) \leq \CC(A)^J \H(F)^J \prod_{i=1}^J H_A(\bx_i),
\end{equation}
where $\CC(A)$ is as in (\ref{aut_const1}). In particular, this implies that
\begin{equation}
\label{prod_4}
H(F X) \leq \H(F)^J \prod_{i=1}^J H(\bx_i).
\end{equation}
\end{lem}

\proof
By Lemmas \ref{wedge} and \ref{comparison},
\begin{equation}
\label{i1}
H_A(F X) = H_{\bigwedge^J A}(\bx_1^t F \wedge\ ...\ \wedge \bx_J^t F) \leq \prod_{i=1}^J H_A(\bx_i^t F) \leq \C_2(A)^J \prod_{i=1}^J H(\bx_i^t F).
\end{equation} 
For each $1 \leq i \leq J$,
$$\bx_i^t F = \left( \sum_{j=1}^N f_{j1} x_{ij}, ..., \sum_{j=1}^N f_{jN} x_{ij} \right).$$
Recall that for the purposes of evaluating height we view the coefficient matrix $F = (f_{ij})_{1 \leq i,j \leq N}$ as a vector in $K^{N^2}$, and we write $\|F\|_v$ for the local norm of this vector at the place $v$. Then for each $v \nmid \infty$,
\begin{equation}
\label{mx_1}
\| \bx_i^t F \|_v \leq \|F\|_v \| \bx_i \|_v,
\end{equation}
and for $v | \infty$, by Cauchy-Schwarz inequality
\begin{eqnarray}
\label{mx_2}
\| \bx_i^t F \|_v & = & \left\{ \sum_{k=1}^N \left\| \sum_{j=1}^N f_{jk} x_{ij} \right\|_v^2 \right\}^{d_v/2d} \nonumber \\
& \leq & \left\{ \sum_{k=1}^N \left( \sum_{j=1}^N \| f_{jk} \|_v^2 \right) \left( \sum_{j=1}^N \| x_{ij} \|_v^2 \right) \right\}^{d_v/2d} = \| F \|_v \| \bx_i \|_v.
\end{eqnarray}
Therefore for each $1 \leq i \leq J$,
\begin{equation}
\label{i2}
H(\bx_i^t F) \leq H(\bx_i) \H(F) \leq \C_1(A)^{-1} H_A(\bx_i) \H(F),
\end{equation}
where the last inequality follows by Lemma \ref{comparison}. Now the lemma follows by combining (\ref{i1}) with (\ref{i2}).
\endproof

\begin{rem} \label{sym_alt} Notice that Lemma \ref{mult} is true for any bilinear form $F$, symmetric, alternating, or none of the above - the proof carries over word for word. Moreover, $F$ can just as well be any $N \times N$ matrix, viewed as a vector in $K^{N^2}$ for the purposes of defining the height $\H(F)$. 
\end{rem}

We are now ready to proceed.

\section{A combinatorial lemma}

In this section we prove a certain graph-theoretic lemma, which we later use in the proof of our main result. We start with some notation. A graph $G$ is connected if there is a path in $G$ connecting every two of its vertices. On the other hand, we will call a pair of vertices connected if they are connected by a single edge, and disconnected otherwise. A graph in which every two vertices are connected is called complete. A complete subgraph on $n$ vertices of a graph $G$ will be called maximal if $G$ does not contain a complete subgraph on any larger number of vertices. Two pairs of vertices in a graph $G$ will be called disjoint if they do not have a vertex in common. We can now state the lemma.

\begin{lem} \label{combin_lemma} Let $G$ be a graph on $2k$ vertices, $k \geq 1$, such that a maximal complete subgraph of $G$ has at most $k$ vertices. Then there exist at least $\left[ \frac{k+1}{2} \right]$ disjoint pairs of disconnected vertices. Moreover, this bound is sharp, meaning that there are such graphs in which any maximal (with respect to cardinality) set of disjoint pairs of disconnected vertices has cardinality precisely $\left[ \frac{k+1}{2} \right]$.
\end{lem}

\proof
Let $v_1,\dots,v_{2k}$ be the vertices of $G$. For each $1 \leq i \neq j \leq 2k$, define
\[ \delta_{ij} = \delta_{ji} = \left\{ \begin{array}{ll}
1 & \mbox{if $v_i$ is connected to $v_j$,} \\
0 & \mbox{otherwise}.
\end{array}
\right. \]
Let 
$$S_1 = \{ 1,\dots,k+1 \},$$
then there must exist $i_1 \neq  j_1 \in S_1$ such that $\delta_{i_1j_1} = 0$: if this was not true, then $G$ would contain a complete subgraph on $k+1$ vertices $v_1,\dots,v_{k+1}$. Next, let
$$S_2 = \left( S_1 \setminus \{ i_1,j_1 \} \right) \cup \{k+2,k+3\}.$$
Since $|S_2|=k+1$, by the same reasoning, there must exist $i_2 \neq  j_2 \in S_2$ such that $\delta_{i_2j_2} = 0$, and next define
$$S_3 = \left( S_2 \setminus \{ i_2,j_2 \} \right) \cup \{k+4,k+5\}.$$
Continuing in this manner, in each set
\begin{equation}
\label{snv}
S_n = \left( S_{n-1} \setminus \{ i_{n-1},j_{n-1} \} \right) \cup \{k+2n-2,k+2n-1\},
\end{equation}
we will find vertices $v_{i_n},v_{j_n}$ such that $\delta_{i_nj_n}=0$. From (\ref{snv}), we see that $1 \leq n \leq M = \left[ \frac{k+1}{2} \right]$, and so we get a collection of distinct vertices
\begin{equation}
\label{vij}
\{ v_{i_1},v_{j_1},\dots,v_{i_M},v_{j_M} \} \subset \{ v_1,\dots,v_{2k} \},
\end{equation}
which satisfy the condition
\begin{equation}
\label{vij1}
\delta_{i_nj_n} = 0,\ \ \forall\ 1 \leq n \leq M = \left[ \frac{k+1}{2} \right].
\end{equation}
This is precisely a collection of $\left[ \frac{k+1}{2} \right]$ disjoint pairs of disconnected vertices in $G$.

Next we show that $\left[ \frac{k+1}{2} \right]$ is sharp. Let $G$ be a graph on vertices $v_1,\dots,v_{2k}$ as above so that $\delta_{ij}=1$ for all $i \neq j$ such that $1 \leq i \leq k-1$ and $1 \leq j \leq 2k$, and $\delta_{ij} = 0$ for all $k \leq i \neq j \leq 2k$; in other words, each of the first $k-1$ vertices is connected to every other vertex in $G$, but no two vertices out of $v_k,\dots,v_{2k}$ are connected to each other. Clearly, any maximal complete subgraph of $G$ will have $k$ vertices; in fact, these will be precisely the $k+1$ subgraphs on the sets of vertices $\{v_1,\dots,v_{k-1},v_j\}$ for each $k \leq j \leq 2k$. Then a maximal (with respect to cardinality) set of disjoint pairs of disconnected vertices is, for instance the set of pairs $v_k,v_{k+1}; \dots; v_{2k-2},v_{2k-1}$, if $k$ is even, and $v_k,v_{k+1}; \dots; v_{2k-1},v_{2k}$, if $k$ is odd. In both cases, the cardinality of such a set is $\left[ \frac{k+1}{2} \right]$. This completes the proof.
\endproof

\section{Proof of Theorem \ref{main}}

In this section we prove a more general version of Theorem \ref{main}, stated as Theorem \ref{main_twisted} below, where the canonical height $H$ is replaced with the twisted height $H_A$, as defined in section~2; since $H$ is simply $H_I$ with $I \in \GL_N(K_{\AA})$ being the identity, Theorem \ref{main} readily follows from Theorem \ref{main_twisted}. We start with a conventional twisted height version of Siegel's lemma.

\begin{thm} \label{siegel} Let $K$ be either a number field, function field, or the algebraic closure of one or the other, and let $Z \subseteq K^N$ be an $L$-dimensional subspace, $1 \leq L < N$. Then for each $A \in \GL_N(K_{\AA})$, there exists a basis $\bz_1,...,\bz_L$ for $Z$ over $K$ such that
\begin{equation}
\label{siegel_bound}
\prod_{i=1}^L H_A(\bz_i) \leq C_K(N,L) H_A(Z),
\end{equation}
where all the notation is as in section~2.

\end{thm}

\proof
When $K$ is a number field, this is the Bombieri-Vaaler version of Siegel's lemma~\cite{vaaler:siegel} with canonical height replaced by twisted height (see \cite{siegel:best}); when $K$ is a function field, this is proved in \cite{thunder}; when $K$ is the algebraic closure of a number field or a function field, this follows from the Roy-Thunder twisted height version of absolute Siegel's lemma (see Theorem 8.1 of \cite{absolute:siegel}).
\endproof

\begin{rem} The constant $C_K(N,L)$ in Theorem \ref{siegel} can be replaced by a slightly sharper one, leading to a slightly better constant in Theorem \ref{main} (see \cite{siegel:best}, \cite{absolute:siegel1}); however, this would make the inequalities harder to read, and some of the constants that would appear in the upper bound would not be easily computable, for instance the generalized Hermite's constant.
\end{rem}

\begin{thm} \label{main_twisted} Let $(Z,F)$ be a regular $2k$-dimensional symplectic space in $N$ variables over $K$, where $1 \leq k < 2k \leq N$. Then for each $A \in \GL_N(K_{\AA})$, there exists a symplectic basis $\bx_1,\dots,\bx_k,\bwy_1,\dots,\bwy_k$ for $Z$ satisfying (\ref{sympl_basis}) such that
\begin{equation}
\label{twisted_basis}
\prod_{i=1}^{2k} H_A(\bx_i)H_A(\bwy_i) \leq \left( C_K(N,2k) H_A(Z) \right)^{a_k}  \left( \CC'(A)  \H(F) \right)^{b_k},
\end{equation}
where $\CC'(A)$ is as in (\ref{aut_const1.1}), and the rest of notation is as in the statement of Theorem \ref{main}. In particular,  if $A = I$ is the identity element of $GL_N(K_{\AA})$, $\CC'(A) = 1$.
\end{thm}

\proof
Fix $A \in \GL_N(K_{\AA})$, and let $\bz_1,\dots,\bz_{2k}$ be the basis for $Z$ guaranteed by Theorem \ref{siegel}. We argue by induction on $k$. If $k=1$, then $F(\bz_1,\bz_2) \neq 0$, since otherwise $(Z,F)$ would be singular. Let $\bx_1 = \frac{1}{F(\bz_1,\bz_2)} \bz_1$, $\bwy_1 = \bz_2$, then $F(\bx_1,\bwy_1) = 1$, and $H_A(\bx_1)=H_A(\bz_1)$. The result follows from (\ref{siegel_bound}). 

Now assume $k>1$. We construct a graph $G(Z)$ on $2k$ vertices in the following way: for each $1 \leq i \leq 2k$, a vertex $v_i$ will correspond to the vector $\bz_i$, and two vertices $v_i$ and $v_j$ will be connected if and only if $F(\bz_i,\bz_j)=0$. Since a Lagrangian of $(Z,F)$ has dimension $k$, the corresponding graph $G(Z)$ satisfies the condition of Lemma \ref{combin_lemma}, which implies that there exists a collection of distinct vectors
\begin{equation}
\label{zij}
\{ \bz_{i_1},\bz_{j_1},\dots,\bz_{i_M},\bz_{j_M} \} \subset \{ \bz_1,\dots,\bz_{2k} \},
\end{equation}
where $M = \left[ \frac{k+1}{2} \right]$, which satisfy the condition
\begin{equation}
\label{zij_nonzero}
F(\bz_{i_n},\bz_{j_n}) \neq 0,\ \ \forall\ 1 \leq n \leq M = \left[ \frac{k+1}{2} \right].
\end{equation}
We can assume without loss of generality that the ordering in (\ref{zij}) satisfies the condition
\begin{equation}
\label{zij_height}
H_A(\bz_{i_1})H_A(\bz_{j_1}) \leq \dots \leq H_A(\bz_{i_M})H_A(\bz_{j_M}).
\end{equation}
Then, combining (\ref{zij_height}) and (\ref{siegel_bound}) we have:
\begin{equation}
\label{zij_height1}
\left( H_A(\bz_{i_1})H_A(\bz_{j_1}) \right)^M \leq \prod_{n=1}^M H_A(\bz_{i_n})H_A(\bz_{j_n}) \leq \prod_{m=1}^{2k} H_A(\bz_m) \leq C_K(N,2k) H_A(Z).
\end{equation}
Let $\bx_1 = \frac{1}{F(\bz_{i_1},\bz_{j_1})} \bz_{i_1}$, $\bwy_1 = \bz_{j_1}$, then $F(\bx_1,\bwy_1) = 1$ and
\begin{equation}
\label{zij_height2}
H_A(\bx_1)H_A(\bwy_1) \leq \left( C_K(N,2k) H_A(Z) \right)^{1/M},
\end{equation}
where $M = \left[ \frac{k+1}{2} \right]$. Let 
$$Z_1 = \spn_K \{ \bx_1, \bwy_1 \}^{\perp_F} \cap Z = \left\{ \bz \in \kbar^N : (\bx_1\ \bwy_1)^t F \bz = \bo \right\} \cap Z,$$ 
then combining Lemmas \ref{dual}, \ref{intersection}, and \ref{mult} with (\ref{aut_const2}), (\ref{comp2}), and (\ref{zij_height2}), we obtain:
\begin{eqnarray}
\label{zij_height3}
H_A(Z_1) & \leq & |\det A|_{\AA} H_{A^*} \left( (\bx_1\ \bwy_1)^t F \right) H_A(Z) \nonumber \\
& \leq & |\det A|_{\AA} \CC(A^*)^2 H_{A^*}(\bx_1) H_{A^*}(\bwy_1) \H(F)^2 H_A(Z) \nonumber \\
& \leq & C_K(N,2k)^{\frac{1}{M}} \left( \frac{\CC(A) |\det A|^{1/2}_{\AA}}{\C_1(A)^2} \right)^2 H_A(Z)^{\frac{M+1}{M}} \H(F)^2.
\end{eqnarray}
Moreover, notice that $\dim_K Z_1 = 2(k-1)$ and $Z_1$ is non-singular, since $Z$ and $\spn_K \{ \bx_1, \bwy_1 \}$ are non-singular. By induction hypothesis, there exists a symplectic basis $\bx_2,\dots,\bx_k,\bwy_2,\dots,\bwy_k$ for $Z_1$ so that
$$F(\bx_i,\bx_j) = F(\bwy_i, \bwy_j) = F(\bx_i,\bwy_j) = 0\ \forall\ 2 \leq i \neq j \leq k,\ F(\bx_i,\bwy_i) = 1\ \forall\ 2 \leq i \leq k,$$
and
\begin{equation}
\label{induct_hyp}
\prod_{i=2}^k H_A(\bx_i)H_A(\bwy_i) \leq \left( C_K(N,2(k-1)) H_A(Z_1) \right)^{a_{k-1}}  \left( \CC'(A) \H(F) \right)^{b_{k-1}}.
\end{equation}
Combining (\ref{zij_height2}), (\ref{zij_height3}), and (\ref{induct_hyp}), and using the fact that $C_K(N,L_1) \leq C_K(N,L_2)$ whenever $L_1 \leq L_2$, we obtain:
\begin{equation}
\label{bnd1}
\prod_{i=1}^k H_A(\bx_i)H_A(\bwy_i) \leq  \left( C_K(N,2k)  H_A(Z) \right)^{\frac{(M+1)a_{k-1}+1}{M}} \left( \CC'(A) \H(F) \right)^{b_{k-1}+2a_{k-1}}.
\end{equation}
The result now follows by a routine calculation.
\endproof

\begin{rem} Clearly, versions of Corollary \ref{hyp_dec} and Corollary \ref{flags} with the twisted height $H_A$ instead of the canonical height $H$ follow immediately from Theorem \ref{main_twisted}.
\end{rem}

\bibliographystyle{plain}  
\bibliography{sympl}        

\end{document}